\documentclass[11pt,reqno]{amsart}
\usepackage{amsmath,amssymb,amsfonts,enumerate,hyperref}
\usepackage{graphicx}
\theoremstyle{plain}
\newtheorem{theorem}{Theorem}[section]
\newtheorem{proposition}[theorem]{Proposition}

\newtheorem{lemma}[theorem]{Lemma}

\theoremstyle{definition}
\newtheorem{definition}[theorem]{Definition}

\theoremstyle{remark}

\numberwithin{equation}{section}

\begin{document}

\title{Quasi-Einstein metrics on hypersurface families}
\author{Stuart James Hall}
\maketitle

\begin{abstract}
We construct quasi-Einstein metrics on some hypersurface families. The hypersurfaces are circle bundles over the product of Fano, K\"ahler-Einstein manifolds.  The quasi-Einstein metrics are related to various gradient K\"ahler-Ricci solitons constructed by Dancer and Wang and some Hermitian, non-K\"ahler, Einstein metrics constructed by Wang and Wang on the same manifolds.
\end{abstract}

\section{Introduction}
\subsection{Motivation and definitions}

This article is concerned with a generalisation of Einstein metrics that in some sense interpolates between Einstein metrics and Ricci solitons, namely, quasi-Einstein metrics. 
\begin{definition}
Let $M^{n}$ be a smooth manifold and $g$ be a complete Riemannian metric. The metric $g$ is called quasi-Einstein if it solves
\begin{equation}\label{qemeq}
Ric(g)+Hess(u)-\frac{1}{m}du \otimes du + \frac{\epsilon}{2} g=0,
\end{equation}
where $u \in C^{\infty}(M)$, $m\in (1,\infty]$ and $\epsilon$ is a constant. 
\end{definition}
It is clear that if $u$ is constant then we recover the notion of an Einstein metric; we will refer to these metrics as trivial quasi-Einstein metrics.  By letting the constant $m$ go to infinity we can also recover the definition of a gradient Ricci soliton. In line with the terminology used for Ricci solitons, we will refer to the quasi-Einstein metrics with $\epsilon< 0 $, $\epsilon = 0$ and $\epsilon >0$ as shrinking, steady and expanding respectively.\\ 

There has been a great deal of effort invested in finding non-trivial examples of Ricci solitons on compact manifolds. However, they remain rare and the only known examples are K\"ahler. Due to work the work of Hamilton \cite{Ham} and Perelman \cite{Per}, non-trivial Ricci solitons on compact manifolds must be shrinking gradient Ricci solitons. The first non-trivial examples were constructed independently by Koiso  and Cao on $\mathbb{CP}^{1}$-bundles over complex projective spaces \cite{Cao, Koi}. These examples were subsequently generalised by Chave and Valet \cite{CV} and Pedersen, T{\o}nneson-Freidman and Valent \cite{PTFV} who found K\"ahler-Ricci solitons on the projectivisation of various line bundles over a Fano K\"ahler-Einstein base. The reader should note that what we call a Ricci soliton is referred to as a quasi-Einstein metrics in the papers \cite{CV} and \cite{PTFV}. Recently Dancer and Wang  generalised these examples by constructing some K\"ahler Ricci solitons on various hypersurface families where the hypersurface is a circle bundle over the product of Fano Kahler-Einstein manifolds \cite{DW1}. The solitons found by Dancer and Wang were also independently constructed by Apostolev, Calderbank, Gauduchon and T{\o}nneson-Freidman \cite{ACGTF}.\\

In the complete non-compact case Feldman, Ilmanen and Knopf \cite{FIK} found shrinking gradient K\"ahler-Ricci solitons on certain line bundles over $\mathbb{CP}^{n}$. Steady gradient K\"ahler-Ricci solitons were first constructed on $\mathbb{C}^{n}$ by Cao \cite{Cao} (the $n=1$ case was first found by Hamilton \cite{Ham2}). Cao also found steady gradient K\"ahler-Ricci solitons on the blow up of $\mathbb{C}^{n}/\mathbb{Z}_{n}$ at the origin. Expanding gradient K\"ahler-Ricci solitons have been found by Cao on  $\mathbb{C}^{n}$ \cite{Cao2} and by Feldman, Ilmanen and Knopf on the blow ups of $\mathbb{C}^{n}/\mathbb{Z}_{k}$ for $k=n+1, n+2,\dots$, \cite{FIK}. Examples were also found by  Pedersen, T\o nneson-Freidman and Valent on the total space of holomorphic line bundles over Kahler-Einstein manifolds with negative scalar curvature \cite{PTFV}. As in the compact case, these examples have been generalised by Dancer and Wang who constructed shrinking, steady and expanding K\"ahler-Ricci solitons on various vector bundles over the product of K\"ahler-Einstein manifolds \cite{DW1}.\\
  
In the recent work \cite{Case10} Case suggested that there should be quasi-Einstein analogues of Dancer-Wang's solitons. He points out that the quasi-Einstein analogue of Koiso-Cao, Chave-Valent and  Pedersen-T\o nneson-Freidman-Valent type solitons was already constructed by L\"u, Page and Pope \cite{LPP}.  The purpose of this article is to show that Dancer-Wang's solitons indeed have quasi-Einstein analogues.  However it is better to think of these metrics as quasi-Einstein analogues of various Hermitian, non-K\"ahler, Einstein metrics constructed by Wang and Wang on these spaces \cite{WW}. The Wang-Wang Einstein metrics generalise a construction originating with Page \cite{Pa} and Berard-Bergery \cite{BB}. We now state the precise results we wish to prove. Non-trivial steady or expanding quasi-Einstein metrics can only occur on non-compact manifolds \cite{KK}.   In the non-compact case we have the following which is the quasi-Einstein analogue of  theorem 1.6 in \cite{WW}:

\begin{theorem}\label{noncom}
Let $(V_{i},J_{i},h_{i}), 1\leq i \leq r, r \geq 3$, be Fano K\"ahler-Einstein manifolds with complex dimension $n_{i}$ and first Chern class $p_{i}a_{i}$ where $p_{i}>0$ and $a_{i}$ are indivisible classes in $H^{2}(V_{i},\mathbb{Z})$. Let $V_{1}$ be a complex projective space with normalised Fubini-Study metric i.e. $p_{1}=(n_{1}+1)$. Let $P_{q}$ denote the principal $\mathbb{S}^{1}$-bundle over $V_{1} \times...\times V_{r}$ with Euler class $\pm \pi_{1}^{\ast}(a_{1})+\sum_{i=2}^{i=r}q_{i}\pi_{i}^{\ast}(a_{i}) $, i.e. $q_{1}^{2}=1$.
\begin{enumerate}

\item Suppose $(n_{1}+1)|q_{i}| <p_{i}$ for $2 \leq i \leq r$ then, for all $m>1$, there exists a non-trivial, complete, steady, quasi-Einstein metric on the total space of the $\mathbb{C}^{n_{1}+1}$-bundle over $V_{2}\times ...\times V_{r}$ corresponding to $P_{q}$.

\item For all $m>1$ there exists at least one one-parameter family of  non-trivial, complete, expanding, quasi-Einstein metrics on the total space of the $\mathbb{C}^{n_{1}+1}$-bundle over $V_{2}\times ...\times V_{r}$ corresponding to $P_{q}$.
\end{enumerate}

\end{theorem}
  
For the compact case we have the following analogue of theorem 1.2 in \cite{WW}.

\begin{theorem}\label{com}
Let $(V_{i},J_{i},h_{i}), 1\leq i \leq r, r \geq 3$, be Fano K\"ahler-Einstein manifolds with complex dimension $n_{i}$ and first Chern class $p_{i}a_{i}$ where $p_{i}>0$ and $a_{i}$ are indivisible classes in $H^{2}(V_{i},\mathbb{Z})$. Let $V_{1}$ and $V_{r}$ be a complex projective space with normalised Fubini-Study metrics. Let $P_{q}$ denote the principal $\mathbb{S}^{1}$-bundle over $V_{1} \times...\times V_{r}$ with Euler class $\pm \pi_{1}^{\ast}(a_{1})+\sum_{i=2}^{i=r-1}q_{i}\pi_{i}^{\ast}(a_{i}) \pm\pi^{\ast}(a_{r})$, i.e. $|q_{1}|=|q_{r}|=1$.\\
 
Suppose that $|q_{i}|(n_{1}+1) < p_{i}$ and  $|q_{i}|(n_{r}+1) < p_{i}$ for $2 \leq i \leq r-1$ and that there exists $\chi = (\chi_{1}, \chi_{2},...,\chi_{r})$ where $|\chi_{i}|=1$, $\chi_{1}=-\chi_{r}=1$ such that
\begin{equation}\label{Inv}
\int_{-(n_{1}+1)}^{(n_{r}+1)}\left(\chi_{1}x+\frac{p_{1}}{|q_{1}|} \right)^{n_{1}}\left(\chi_{2}x+\frac{p_{2}}{|q_{2}|} \right)^{n_{2}}...\left(\chi_{r}x+\frac{p_{r}}{|q_{r}|} \right)^{n_{r}}xdx < 0,
\end{equation}
then, for all $m>1$ there exists a non-trivial, shrinking quasi-Einstein metric on $M_{q}$, the space obtained from $P_{q}\times_{S_{1}} \mathbb{CP}^{1}$ by blowing-down one end to $V_{2} \times ...\times V_{r}$ and the other end to $V_{1} \times ...\times V_{r-1}$. 

\end{theorem}

We remark that the Futaki invariant (evaluated on the holomorphic vector field $f(t)\partial_{t}$ in the notation of the next section) is given by
$$ \int_{-(n_{1}+1)}^{(n_{r}+1)}\left(\frac{p_{1}}{q_{1}}-x \right)^{n_{1}}\left(\frac{p_{2}}{q_{2}}-x \right)^{n_{2}}...\left(\frac{p_{r}}{q_{r}}-x \right)^{n_{r}}xdx.$$ If this integral vanishes then Dancer-Wang construct a K\"ahler-Einstein metric on $M_{q}$. 

Finally we note that none of the metrics we find are K\"ahler. Indeed there is a rigidity result due to Case-Shu-Wei \cite{CSW} that says, on compact manifolds, K\"ahler-quasi-Einstein metrics are trivial i.e. K\"ahler-Einstein.\\
\\
\textit{Acknowledgements: } I would like to thank Prof. Andrew Dancer for many interesting conversations about quasi-Einstein metrics and Ricci solitons.  I would also like to thank Maria Buzano, Jeffrey Case and Tommy Murphy for useful comments on this paper. I would also like to thank the anonymous referee for useful suggestions and corrections to the previous version.
\section{Proof of main theorems}
\subsection{Derivation of equations}
We use the same notation as above. We consider the manifold \linebreak $M_{0} = (0,l)\times P_{q}$.  Let $\theta$ be the principal $U(1)$-connection on $P_{q}$ with curvature $\Omega = \sum_{i=1}^{r}q_{i}\pi^{\ast}\eta_{i}$ where $\eta_{i}$ is the K\"ahler form of the metric $h_{i}$. We form the 1-parameter family of metrics on $P_{q}$
$$g_{t} = f^{2}(t)\theta \otimes \theta +\sum_{i=1}^{i=r} g_{i}^{2}(t)\pi^{\ast}h_{i}$$
and we then form the metric $\bar{g} =dt^{2}+g_{t}$ on $M_{0}$. The group $U(1)$ acts on $M_{0}$ by isometries and generates a Killing field $Z$.  We define a complex structure $J$ on $M_{0}$ by $J(\partial_{t})=-f^{-1}(t)Z$ on the vertical space of $\theta$ and by lifting the complex structure from the base on the horizontal spaces of $\theta$.

\begin{lemma}
Let $M_{0}$ be as above and let $v=e^{-\frac{u}{m}}$. Then the quasi-Einstein equations in this setting are given by:
\begin{equation}\label{QE1}
\frac{\ddot{f}}{f}+\sum_{i=1}^{i=r}2n_{i}\frac{\ddot{g}_{i}}{g_{i}}+m\frac{\ddot{v}}{v}  =  \frac{\epsilon}{2},
\end{equation}
\begin{equation}\label{QE2}
\frac{\ddot{f}}{f}+\sum_{i=1}^{i=r}\left(2n_{i}\frac{\dot{f}\dot{g}_{i}}{fg_{i}}-\frac{n_{i}q^{2}_{i}}{2}\frac{f^{2}}{g_{i}^{4}} \right)+ m\frac{\dot{f}\dot{v}}{fv}  =  \frac{\epsilon}{2},
\end{equation}
\begin{equation}\label{QE3}
\frac{\ddot{g}_{i}}{g_{i}}-\left( \frac{\dot{g}_{i}}{g_{i}}\right)^{2}+\frac{\dot{f}\dot{g}_{i}}{fg_{i}}+\sum_{j=1}^{j=r}2n_{j}\frac{\dot{g}_{i}\dot{g}_{j}}{g_{i}g_{j}}-\frac{p_{i}}{g_{i}^{2}}+\frac{q_{i}^{2}f^{2}}{2g_{i}^{4}}+ m\frac{\dot{g_{i}}\dot{v}}{g_{i}v} = \frac{\epsilon}{2}.
\end{equation}
\end{lemma}

In order that $(M,g,u)$ be a quasi-Einstein manifold, as well as equation (\ref{qemeq}), $u$ must also satisfy an integrability condition that essentially comes from the second Bianchi identity (c.f. Lemma 3.4 in \cite{DW1}). The form we use here is given in Case \cite{Case10}  using the Bakry-\'Emery Laplacian:
$$\Delta_{u}:=\Delta - \langle\nabla u, \cdot \rangle.$$
\begin{lemma}[Kim-Kim \cite{KK} Corollary 3]
Let $(M,g,u)$ be a quasi-Einstein manifold then there exists a constant $\mu$ such that
\begin{equation}\label{mueq}
\Delta_{u}\left(\frac{u}{m}\right)+\frac{\epsilon}{2}=-\mu e^{\frac{2u}{m}}.
\end{equation}
In the notation above (recalling  $v=e^{-\frac{u}{m}}$) this condition becomes 
\begin{equation}\label{int1}
\mu = v\ddot{v}+v\dot{v}\left(\frac{\dot{f}}{f}+\sum_{i}2n_{i}\frac{\dot{g}_{i}}{g_{i}}\right)+(m-1)\dot{v}^{2}-\frac{\epsilon}{2}v^{2}.
\end{equation}
\end{lemma}
The constant $\mu$ enters into the discussion of Einstein warped products when $m$ is an integer.  If $(M,g,u)$ is a quasi-Einstein manifold with constant $\mu$ coming from (\ref{mueq}) and $(F^{m},h)$ is an Einstein manifold with constant $\mu$, then $(M\times F^{m},g\oplus v^{2}h)$ is an Einstein metric with constant $-\epsilon/2$ as in equation (\ref{qemeq}) (c.f. \cite{KK}).\\
\\
Introducing the moment map change of variables as in \cite{DW1} and \cite{WW} yields the following set of equations:
\begin{proposition}
Let $s$ be the coordinate on $I=(0,l)$ such that $ds=f(t)dt$, $\alpha(s)=f^{2}(t)$, $\beta_{i}(s)=g_{i}^{2}(t)$, $\phi(s)=v(t)$ and $V = \prod_{i=1}^{i=r}g_{i}^{2n_{i}}(t).$ Then the equations (\ref{QE1}),(\ref{QE2}),(\ref{QE3}) and (\ref{int1}) transform to the following:
\begin{equation}\label{mmeq1}
\frac{1}{2}\alpha'' +\frac{1}{2}\alpha'(\log V)' +\alpha\sum_{i=1}^{r}n_{i}\left(\frac{\beta_{i}''}{\beta_{i}}-\frac{1}{2}\left(\frac{\beta_{i}'}{\beta_{i}}\right)^{2}\right)+m\left( \frac{\alpha \phi''}{\phi}+\frac{\alpha'\phi'}{2\phi} \right)=\frac{\epsilon}{2},
\end{equation}
\begin{equation}\label{mmeq2}
\frac{1}{2}\alpha'' +\frac{1}{2}\alpha'(\log V)'-\alpha\sum_{i=1}^{i=r}\frac{n_{i}q_{i}^{2}}{2\beta_{i}^{2}}+m\frac{\alpha'\phi'}{2\phi}=\frac{\epsilon}{2},
\end{equation}
\begin{equation}\label{mmeq3}
\frac{1}{2}\frac{\alpha'\beta_{i}'}{\beta_{i}}+\frac{1}{2}\alpha\left(\frac{\beta_{i}''}{\beta_{i}}-\left(\frac{\beta_{i}'}{\beta_{i}} \right)^{2} \right)+\frac{1}{2}\frac{\alpha \beta_{i}'}{\beta_{i}}(\log V)' -\frac{p_{i}}{\beta_{i}}+\frac{q_{i}^{2}\alpha}{2\beta_{i}^{2}}+m\frac{\alpha}{2}\frac{\beta_{i}'\phi'}{\beta_{i}\phi}=\frac{\epsilon}{2},
\end{equation}
\begin{equation}\label{inteq}
\phi\left(\phi''\alpha +\frac{\phi'\alpha'}{2}\right)+\phi\phi'\left( \frac{\alpha'}{2}+ (\log V)'\alpha  \right)+(m-1)(\phi')^{2}\alpha-\frac{\epsilon}{2}\phi^{2}=\mu.
\end{equation}
\end{proposition}
 Equating (\ref{mmeq1}) and (\ref{mmeq2}) we obtain
\begin{equation}
-m\frac{\phi''}{\phi} = \sum_{i=1}^{i=r}n_{i}\left(\frac{\beta_{i}''}{\beta_{i}}-\frac{1}{2}\left(\frac{\beta_{i}'}{\beta_{i}} \right)^{2}+\frac{q_{i}^{2}}{2\beta_{i}^{2}} \right)
\end{equation}
Following \cite{DW1,WW} we look for solutions that satisfy
$$\frac{\beta_{i}''}{\beta_{i}}-\frac{1}{2}\left(\frac{\beta_{i}'}{\beta_{i}}\right)^{2}+\frac{1}{2}\frac{q_{i}^{2}}{\beta_{i}^{2}}=0.$$
This condition can be geometrically interpreted as saying that the curvature of $\overline{g}$ is $J$-invariant in the sense that $\overline{Rm}(J\cdot,J\cdot,J\cdot,J\cdot) = \overline{Rm}(\cdot,\cdot,\cdot,\cdot)$ where $J$ is the complex structure on $M_{0}$. Imposing this forces $\phi$ to be a linear function of $s$.  We write $\phi(s) = \kappa_{1}(s+\kappa_{0})$ for some constants $\kappa_{0},\kappa_{1} \in \mathbb{R}$. Hence (\ref{inteq}) becomes
\begin{equation}\label{inteq2}
\alpha '+\alpha((\log V)' +\frac{(m-1)}{(s+\kappa_{0})}) =\frac{\epsilon (s+\kappa_{0})}{2}+\frac{\mu}{\kappa_{1}^{2}(s+\kappa_{0})}.
\end{equation}
Accordingly there are two classes of solution $\beta_{i}$: 
$$\beta_{i}=A_{i}(s+s_{0})^{2}-\frac{q_{i}^{2}}{4A_{i}}$$ or
$$\beta_{i} = \pm q_{i}(s+\sigma_{i})$$
where $A_{i} \neq 0$ and $\sigma_{i}$ are constants.  We note that the case $\beta_{i}=-q_{i}(s+\sigma_{i})$ corresponds to the metric $\bar{g}$ being K\"ahler with respect to the complex structure. Hence the rigidity result of Case-Shu-Wei rules out having any solutions of this form (in fact choosing $\beta_{i}$ of this form leads to inconsistency).\\

If we input $\beta_{i} = A_{i}(s+s_{0})^{2}-\frac{q_{i}^{2}}{4A_{i}}$ into (\ref{mmeq3}) we see that
$$\alpha'+\alpha\left((\log V)'+m (\log \phi)'-\frac{1}{(s+s_{0})}\right) =  \frac{\epsilon}{2}(s+\kappa_{0})+\frac{E^{\ast}}{(s+\kappa_{0})}$$
where $$E^{\ast} :=\frac{8A_{i}p_{i}-\epsilon q_{i}^{2}}{8A_{i}^{2}}.$$ 
Comparing with equation (\ref{inteq2}) we see that solutions are consistent providing
$\kappa_{0}=s_{0}$  and 
$$ \frac{\mu}{\kappa^{2}_{1}}=E^{\ast}=\frac{8A_{i}p_{i}-\epsilon q_{i}^{2}}{8A_{i}^{2}} .$$
Solving gives
\begin{equation}
\alpha(s) = V^{-1}(s+\kappa_{0})^{1-m}\int_{0}^{s}V(s+\kappa_{0})^{m-2}\left(E^{\ast}+\frac{\epsilon}{2}(s+\kappa_{0})^{2}\right)ds.
\end{equation}

\subsection{Compactifying $M_{0}$}
We recall that $V_{1} = \mathbb{CP}^{n_{1}}$ and we are adding in the manifold $V_{2}\times...\times V_{r}$ at the point $s=0$.  We refer the reader to the discussion immediately after equation (4.17)  in \cite{DW1}.  In a nutshell, in order for the metric to extend smoothly we require that 
$$\alpha(0)=0,\alpha'(0)=2, \beta_{1}(0)=0 \textrm{ and } \beta_{1}'(0)=1.$$
As we are using $\beta_{1}(s) = A_{1}(s+\kappa_{0})^{2}-\frac{q_{1}^{2}}{4A_{1}}$ we must have $A_{1} = \frac{1}{2\kappa_{0}}$ and $|q_{1}|=1$. We also have normalised so that $p_{1}=n_{1}+1$ hence the consistency conditions become
$$E^{\ast} = \frac{\mu}{\kappa_{1}^{2}} = \frac{\kappa_{0}}{2}(4(n_{1}+1)-\epsilon\kappa_{0})=\frac{8A_{i}p_{i}-\epsilon q_{i}^{2}}{8A_{i}^{2}} \textrm{ for } 2 \leq i \leq r.$$
\subsection{Steady quasi-Einstein metrics}
In this case $\epsilon =0$. Setting $V_{1} =\mathbb{CP}^{n_{r}}$ and compactifying we obtain a $\mathbb{C}^{n_{1}+1}$-vector bundle over $V_{2}\times ...\times V_{r}$.  In order that $\beta_{i}(0)>0$ on $I=[0,\infty)$ we must have $A_{i}>0$ and
$$E^{\ast} = \frac{\mu}{\kappa_1^{2}} = \frac{\kappa_{0}}{2}(4(n_{1}+1))=\frac{p_{i}}{A_{i}} \textrm{ for } 2 \leq i \leq r.$$ Hence $A_{i} = \frac{p_{i}}{E^{\ast}}$ and
$$\beta_{i}(s) = \frac{p_{i}}{E^{\ast}}(s+\kappa_{0})^{2}-\frac{E^{\ast}q_{i}^{2}}{4p_{i}}.$$
It is clear that in order for $\beta_{i}(0)>0$ we must have
$$(n_{1}+1)|q_{i}| <p_{i}  \textrm{ for } 2 \leq i \leq r.$$ In order to ensure the metrics are complete it is sufficient to check that the integral
\begin{equation}\label{compint}
t = \int_{0}^{s}\frac{dx}{\sqrt{\alpha(x)}}
\end{equation}
diverges as $s \rightarrow \infty$ (this says that geodesics cannot reach the boundary at infinity and have finite length).  If we compute the function $\alpha(s)$ we see that it is asymptotic to a positive constant $K$. Hence the above integral diverges and we obtain a complete quasi-Einstein metric for all $m>1$ generalising the non-K\"ahler, Ricci-flat ones constructed in \cite{WW}. Choosing a different value of $E^{\ast}$ simply varies the metric by homothety. 
\subsection{Expanding quasi-Einstein metrics}
Here we take $\epsilon=1$ to factor out homothety.  Again the manifolds in question are $\mathbb{C}^{n_{1}+1}$-vector bundles over $V_{2}\times ...\times V_{r}$.
Here the consistency conditions become
$$E^{\ast} = \frac{\mu}{\kappa_1^{2}} = \frac{\kappa_{0}}{2}(4(n_{1}+1)-\kappa_{0})=\frac{8A_{i}p_{i}- q_{i}^{2}}{8A_{i}^{2}} \textrm{ for } 2 \leq i \leq r.$$
 If $|q_{i}|(n_{1}+1) \leq p_{i}$ then we choose $0<E^{\ast}< 2(n_{1}+1)^{2}$,
 $$\kappa_{0} = 2(n_{1}+1)+2\sqrt{(n_{1}+1)^{2}-\frac{E^{\ast}}{2}} $$
 and
 $$A_{i} = \frac{1}{2E^{\ast}}\left(p_{i}+\sqrt{p_{i}^{2}-\frac{E^{\ast}q_{i}^{2}}{2}}\right).$$
In order that $\beta_{i}(0)>0$ we require $2\kappa_{0}A_{i}>|q_{i}|$ for $2 \leq i\leq r$. This can be seen as
$$2\left( 2(n_{1}+1)+2\sqrt{(n_{1}+1)^{2}-\frac{E^{\ast}}{2}}\right)\frac{1}{2E^{\ast}}\left(p_{i}+\sqrt{p_{i}^{2}-\frac{E^{\ast}q_{i}^{2}}{2}}\right)>\frac{2(n_{1}+1)p_{i}}{E^{\ast}}>|q_{i}|.$$
In the case that $|q_{i}|(n_{1}+1)<p_{i}$ we note also that
$$\left(1+\sqrt{1-\frac{E^{\ast}q_{i}^{2}}{2p_{i}^{2}}}\right)>\left( 1+\sqrt{1-\frac{E^{\ast}}{2(n_{1}+1)^{2}}}\right),$$
hence,
$$ 2\left( 2(n_{1}+1)-2\sqrt{(n_{1}+1)^{2}-\frac{E^{\ast}}{2}}\right)\frac{1}{2E^{\ast}}\left(p_{i}+\sqrt{p_{i}^{2}-\frac{E^{\ast}q_{i}^{2}}{2}}\right) > $$
$$\frac{4p_{i}(n_{1}+1)}{2E^{\ast}}\left(1-\sqrt{1-\frac{E^{\ast}}{2(n_{1}+1)^{2}}} \right)\left(1+\sqrt{1+\frac{E^{\ast}}{2(n_{1}+1)^{2}}} \right)= \frac{p_{i}}{(n_{1}+1)}>|q_{i}|.$$
Therefore if we have the strict inequality $|q_{i}|(n_{1}+1)<p_{i}$ then we can also choose
$$\kappa_{0} = 2(n_{1}+1)-2\sqrt{(n_{1}+1)^{2}-\frac{E^{\ast}}{2}}.$$
If  $|q_{i}|(n_{1}+1)>p_{i}$ then we can choose $0<E^{\ast} < 2(n_{1}+1)^{2}\min (p_{2}^{2}/q_{2}^{2},...,p_{r}^{2}/q_{r}^{2})$.  If we also choose
 $$\kappa_{0} = 2(n_{1}+1)+2\sqrt{(n_{1}+1)^{2}-\frac{E^{\ast}}{2}}$$ and 
 $$A_{i} = \frac{1}{2E^{\ast}}p_{i}+\sqrt{p_{i}^{2}-\frac{E^{\ast}q_{i}^{2}}{2}},$$
 then $\beta_{i}(0)>0$. We can also choose $E^{\ast}<0$ in this case.  Completeness follows as $\alpha(s)$ is asymptotic to $Ks^{2}$ for a positive constant $K$ and so the integral (\ref{compint}) diverges. Hence we find complete, quasi-Einstein analogues of the non-Kahler, Einstein metrics constructed in \cite{WW}. 
\subsection{Shrinking quasi-Einstein metrics}
In order to factor out homothety we take $\epsilon=-1$ and so the consistency conditions are
$$\frac{\mu}{\kappa_{1}^{2}} = \frac{\kappa_{0}}{2}(4(n_{1}+1)+\kappa_{0})=\frac{8A_{i}p_{i}+ q_{i}^{2}}{8A_{i}^{2}} \textrm{ for } 2 \leq i \leq r.$$
We split the discussion into the compact case and the non-compact, complete case.  For the compact case we consider $I$ to be the finite interval $[0, s_{\ast}]$. We set $V_{r} = \mathbb{CP}^{n_{r}}$ and at the point $s=s_{\ast}$ we add in the manifold $V_{1}\times...\times V_{r-1}$. For the metric to extend smoothly we require that $q_{r}=1, p_{r}=n_{r}+1$ and $-1=2A_{r}(s_{\ast}+\kappa_{0})$. Putting these into the consistency conditions we see that
$$ \kappa_{0}(4(n_{1}+1)+\kappa_{0}) = (s_{\ast}+\kappa_{0})^{2}-4(n_{r}+1)(s_{\ast}+\kappa_{0})$$
and hence
$$s_{\ast} =  \sqrt{\kappa_{0}(4(n_{1}+1)+\kappa_{0})+4(n_{r}+1)^{2}}-\kappa_0+2(n_{r}+1).$$
We note that if $n_{1}=n_{r}$ then $s_{\ast} = 4(n_{1}+1)$.
For the time being we note that $s_{\ast}=s_{\ast}(E^{\ast})$ and $\beta_{i}$ is completely determined by $E^{\ast}$ once we have chosen the value of $q_{i}^{2}$ and the sign of $A_{i}$.
The value $A_{i}$ is given by
$$A_{i} = \frac{1}{2E^{\ast}}\left(p_{i} + \chi_{i}  \sqrt{p_{i}^{2}+\frac{E^{\ast}q_{i}^{2}}{2}}\right)$$
where $\chi_{i} =1$ if $A_{i}>0$ and $\chi_{i} =-1$ if $A_{i}<0$.
In order to have a quasi-Einstein metric we must be able  choose a value  of $E^{\ast}>0$ such that the integral
$$\int_{0}^{s_{\ast}(E^{\ast})}\prod_{i=0}^{i=r}\left[\left((s+\kappa_{0})^{2}-\frac{q_{i}^{2}}{4A_{i}^{2}}\right)^{n_{i}}\right](s+\kappa_{0})^{m-2}\left(E^{\ast}-\frac{1}{2}(s+\kappa_{0})^{2}\right)ds = 0.$$
Changing coordinates to
$$x = \frac{1}{2}(s+\kappa_0) - ((n_{1}+1)^{2}+\frac{E^{\ast}}{2})^{1/2},$$
then the above integral becomes (ignoring constants)
$$ F(E^{\ast}) = \int_{-(n_{1}+1)}^{x_{\ast}(E^{\ast})}\prod_{i=0}^{i=r}P_{i}(x)(x+((n_{1}+1)^{2}+\frac{E^{\ast}}{2})^{1/2})^{m-2}(x^{2}+2x((n_{1}+1)^{2}+\frac{E^{\ast}}{2})^{1/2})+(n_{1}+1)^{2})ds$$
where
$$P_{i}(x) = \left(x^{2}+ 2x((n_{1}+1)^{2}+\frac{E^{\ast}}{2})^{1/2}+(n_{1}+1)^{2}+\frac{2p_{i}(\chi_{i}\sqrt{p_{i}^{2}+\frac{E^{\ast}q_{i}^{2}}{2}}-p_{i})}{q_{i}^{2}}\right)^{n_{i}}$$
and
$$x_{\ast}(E^{\ast}) =(n_{r}+1)+(\frac{E^{\ast}}{2}+(n_{r}+1)^{2})^{1/2}-(\frac{E^{\ast}}{2}+(n_{1}+1)^{2})^{1/2}.$$
We will compute the limit $\lim_{E^{\ast}\downarrow 0} F(E^{\ast})$ and the limit $\lim_{E^{\ast} \rightarrow \infty} F(E^{\ast})$.  

We begin with $0$. We note that as $m>1$ the function
$f(x) = (x+(n_{1}+1)^{m-2}$ is integrable on $[-(n_{1}+1),x(E^{\ast})]$ so by the dominated convergence theorem we can evaluate the integral of the limit. This is given by
$$S\int_{-(n_{1}+1)}^{2(n_{r}+1)-(n_{1}+1)}\prod_{\chi_{i}=-1}\left[x+(n_{1}+1)\right]^{2n_{i}}\prod_{\chi_{j}=1}\left[\frac{4p_{i}^{2}}{q_{i}^{2}} -(x+(n_{1}+1))^{2}\right]^{n_{j}}(x+(n_{1}+1))^{m} dx,$$ 
where 
$$S = (-1)^{\sum_{\chi_{i}=-1}n_{i}}.$$
The hypothesis on the $p_{i}$ and $q_{i}$ mean that the sign of $\lim_{E^{\ast}\downarrow 0} F(E^{\ast})$ is that of $S$.\\

For $E^{\ast} \rightarrow \infty$ we consider 
$$\lim_{E^{\ast} \rightarrow \infty}F(E^{\ast})(E^{\ast})^{\frac{1}{2}(1-m-\sum_{\chi_{i}=-1}n_{i})} =K(-1)^{\sum_{\chi_{i}=-1}n_{i}}\int_{-(n_{1}+1)}^{(n_{r}+1)}\prod_{i=1}^{i=r}\left[\chi_{i}x+\frac{p_{i}}{|q_{i}|}\right]^{n_{i}}xdx,$$
 where $K$ is a positive constant. Hence if we can choose $\chi_{i}$ so that 
$$ \int_{-(n_{1}+1)}^{(n_{r}+1)}\prod_{i=1}^{i=r}\left[\chi_{i}x+\frac{p_{i}}{|q_{i}|}\right]^{n_{i}}xdx<0,$$ we can find an $E^{\ast}>0$ such that $\alpha(s_{\ast})=0$.  A discussion similar to that in \cite{DW1} and \cite{WW} shows that this is enough to ensure we have smooth quasi-Einstein metrics.

\section{Examples and future work}
We end with an example of theorem \ref{com}, some discussion of the geometry of the quasi-Einstein metrics constructed and a discussion of possible sources future compact examples.
\subsection{An example}
We consider an example that is also considered in \cite{DW1}.  They consider a $\mathbb{CP}^{1}$-bundle over  $\mathbb{CP}^{2}\times\mathbb{CP}^{2}$.  In theorem \ref{com} this corresponds to taking $r=4, n_{1}=n_{4}=0, n_{2}=n_{3}=2$ and $p_{2}=p_{3}=3$. Hence to apply the theorem we must consider $|q_{2}|,|q_{3}|<3$.  They take $(q_{2},q_{3}) = (1,-2).$
The Futaki invariant is given by
$$\int_{-1}^{1}(3-x)^{2}(\frac{3}{2}+x)^{2}xdx$$
which they calculate is  $7.8$.   This means that
$$\int_{-1}^{1}(3+x)^{2}(\frac{3}{2}-x)^{2}xdx =-7.8<0$$
and we have non-trivial quasi-Einstein metrics on this space for all $m>1$.

\subsection{Remarks on the geometry of the quasi-Einstein metrics}

In \cite{DW1} section 4, the authors comment on the geometry at infinity of their examples of steady and expanding gradient K\"ahler-Ricci solitons.  In particular they conclude that their steady examples are asympotically parabolic and that the expanding examples are asymptotically conical.  We recall that the examples of steady quasi-Einstein metrics constructed in theorem \ref{noncom} have $\alpha(s) \sim K$ for some positive constant $K$ and so the following asymptotic behaviour holds (ignoring multiplicative constants)
$$ f(t) =O(1) \textrm{ and }  g_{i}(t) \sim t.$$
In the expanding case we recall that $\alpha(s) \sim Ks^{2}$ and so we have
$$ f(t) \sim e^{t} \textrm{ and }  g_{i}(t) \sim e^{t}.$$

\subsection{Future families}

The space $\mathbb{CP}^{2} \sharp \overline{\mathbb{CP}}^{2}$ fits into the framework of theorem \ref{com} as a non-trivial $\mathbb{CP}^{1}$-bundle over $\mathbb{CP}^{1}$. On this space  there is the Page metric, the Koiso-Cao soliton and the quasi-Einstein metrics of theorem $3$ (originally due to L\"u-Page-Pope).  The space $\mathbb{CP}^{2} \sharp 2 \overline{  \mathbb{CP}}^{2}$ also admits a non-K\"ahler, Einstein metric due to Chen, LeBrun and Weber \cite{CLB} and a K\"ahler-Ricci soliton due to Wang and Zhu \cite{WZ}.  It would seem reasonable that there should be a family of quasi-Einstein analogues to these metrics.  The metrics on $\mathbb{CP}^{2} \sharp 2 \overline{\mathbb{CP}}^{2}$  are not cohomogeneity-one but do have an isometric action by $\mathbb{T}^{2}$. One observation is that the L\"u-Page-Pope quasi-Einstein metrics are conformally K\"ahler (as any $U(2)$-invariant metric on $\mathbb{CP}^{2}\sharp\overline{\mathbb{CP}}^{2}$ is). The Chen-LeBrun-Weber metric is also conformally K\"ahler (a fact orginally proved by Derdzinski \cite{Der}) and so one might hope that the same would be true for analogues of the L\"u-Page-Pope metrics on  $\mathbb{CP}^{2}\sharp2\overline{\mathbb{CP}}^{2}$.  Both the Page and Chen-LeBrun-Weber metrics are conformal to extremal K\"ahler metrics which satisfy a PDE that `occurs naturally' in K\"ahler geometry.  It would be an interesting first step to try and find an analogous PDE/ODE for the K\"ahler metrics that are conformal to the L\"u-Page-Pope metrics. The author hopes to take up the existence questions in a future work.

\address{Department of Applied Computing, The University of Buckingham, Hunter Street
Buckingham, MK18 1EG,United Kingdom\\
email: stuart.hall@buckingham.ac.uk}

\begin{thebibliography}{999}

\bibitem{ACGTF} V. Apostolev, D. Calderbank, P. Gauduchon, C. T{\o}nnesen-Friedman, Hamiltonian 2-forms in K\"ahler geometry IV: Weakly Bochner-flat K\"ahler manifolds, Commum. Anal. Geom., 16, (2008), 91--126.

\bibitem{BB} L. B\'erard-Bergery, Sur des nouvelles vari\'et\'es riemannienes d'Einstein, Publication de l'Institute Elie Cartan, Nancy (1982).

\bibitem{Cao} H.-D. Cao, Existence of gradient Ricci solitons, Elliptic and parabolic methods in geometry,A.K. Peters, Wellesley, (1996) 1-16.

\bibitem{Cao2} H.-D. Cao, Limits of solutions to the K\"ahler-Ricci flow, J. Diff, Geom., 45, (1997), 257--272.

\bibitem{CSW} J. Case, Y. Shu, G Wei, Rigidity of quasi-Einstein metrics, Diff. Geom Appl., 29, (2011), 93--100.
 
\bibitem{Case10} J. Case, Smooth metric measure spaces and quasi-Einstein metrics, preprint, (2010), arXiv:1011.2723v3 [math.DG].

\bibitem{CV} T. Chave, G. Valent, On a class of compact and non-compact quasi-Einstein metrics and their renormalizability properties, Nuclear Phys. B, 478,  no. 3, (1996), 758--778. 

\bibitem{CLB} X. Chen, C. LeBrun, B. Weber, On conformally K\"ahler, Einstein manifolds, J. Amer. Math. Soc., 21, no. 4, (2008), 1137--1168.

\bibitem{DW1} A. Dancer, M. Wang, On Ricci solitons of cohomogeneity one, Annals of Global Analysis and Geometry, 39, (2011), 259--292.

\bibitem{Der} A. Derdzinski, Self-dual K\"ahler manifolds and Einstein manifolds of dimension four, Compositio Math., 49, no. 3, (1983), 405--433.

\bibitem{FIK} M. Feldman, T. Ilmanen, D. Knopf, Rotationally symmetric shrinking and expanding gradient K\"ahler-Ricci solitons, J. Diff. Geom., 65, (2003), 169--209.

\bibitem{Ham2} R. Hamilton, The Ricci flow on surfaces, Contemporary Mathematics, 71, (1988), 237--261.

\bibitem{Ham} R. Hamilton, The formation of singularities in the Ricci flow, Surveys in differential geometry, Voll. II (Cambridge, MA, 1993), (1995), 7--136.

\bibitem{KK} D.-S. Kim, Y.H. Kim, Compact Einstein warped product spaces with nonpositive scalar curvature, Proc. Amer. Math. Soc. 131, (2003), 2573--2576.

\bibitem{Koi}N. Koiso,  On rotationally symmetric Hamilton's equation for K\"ahler-Einstein metrics, Advanced studies in Pure Mathematics, vol. 18-I. Academic Press, Tokyo, (1990), 327--337.

\bibitem{LPP}H. L\"u, D. Page, C. Pope, New inhomogenous Einstein metrics on sphere bundles over Einstein-K\"ahler manifolds, Phys Lett B, 593, (2004), 218--226.

\bibitem{Pa} D. Page, A compact rotating gravitational instanton, Phys Lett, 79B, (1979), 235--238.

\bibitem{PTFV} H. Pedersen, C. T\o nnesen-Freidman, G. Valent, Quasi-Einstein K\"ahler metrics, Lett. Math. Phys., 50, no. 3, (1999), 229--241.

\bibitem{Per} G. Perelman, The entropy formula for the Ricci flow and its geometric applications, preprint, (2002), arXiv:math/0211159 [math.DG].

\bibitem{WW} J. Wang, M. Wang, Einstein metrics on $S^{2}$-bundles, Math Ann, 310, (1998), 497--526. 

\bibitem{WZ} X.-J.  Wang, X. Zhu, K\"ahler-Ricci solitons on toric manifolds with positive first Chern class. Adv.  Math. 188, (2004), 87--103. 


\end{thebibliography}
\end{document}